\documentclass[11pt]{article}
\usepackage[english]{babel}
\usepackage{newlfont}

\usepackage{amsthm}
\usepackage{hyperref}
\usepackage{mathtools}
\usepackage{amsfonts}
\usepackage{amssymb}
\usepackage{esint}
\usepackage{xcolor}
\usepackage[shortlabels]{enumitem}
\usepackage{bigints}
\usepackage{dsfont}
\usepackage{mathrsfs}
\usepackage{ragged2e}
\usepackage{parskip} 

\usepackage{amsmath}
\usepackage{tikz}
\usepackage{mathdots}
\usepackage{cancel}
\usepackage{color}
\usepackage{siunitx}
\usepackage{array}
\usepackage{multirow}
\usepackage{amssymb}
\usepackage{tabularx}
\usepackage{extarrows}
\usepackage{booktabs}
\usetikzlibrary{fadings}
\usetikzlibrary{patterns}
\usetikzlibrary{shadows.blur}
\usetikzlibrary{shapes}

\usepackage{geometry}
\theoremstyle{definition}
\newtheorem{defn}{Definition}[section]
\newtheorem{thm}[defn]{Theorem}

\newtheorem{rem}[defn]{Remark}
\newtheorem{propx}[defn]{Proposition}
\newtheorem{cor}[defn]{Corollary}
\newtheorem{lemma}[defn]{Lemma}

\numberwithin{equation}{section}

\newcommand{\R}{\mathbb{R}}
\newcommand{\C}{\mathbb{C}}
\newcommand{\N}{\mathbb{N}}
\newcommand{\Z}{\mathbb{Z}}

\newcommand{\Hyp}{\mathbb{H}}

\newcommand{\eps}{\varepsilon}
\newcommand{\la}{\langle}
\newcommand{\ra}{\rangle}

\DeclareMathOperator{\arccosh}{arccosh}
\DeclareMathOperator{\Pl}{Pl}
\DeclareMathOperator{\InjRad}{InjRad}
\DeclareMathOperator{\SL}{SL}

\title{\huge{\textbf{Benjamini--Schramm vs Plancherel convergence}}}
\author{{\LARGE Giacomo Gavelli, Claudius Kamp}}
\date{}

\begin{document}
\maketitle

{\large \textbf{Abstract: }We show that Plancherel convergence is strictly stronger than Benjamini -- Schramm convergence. In order to do so, we introduce two criteria for Plancherel and Benjamini-Schramm convergence in terms of counting functions of the length spectrum.}

\tableofcontents
\newpage

\section{Introduction}
Let $G$ be a connected semisimple Lie group of non-compact type, $K\le G$ a maximal compact subgroup and let $X=G/K$ the associated Riemannian symmetric space. A lattice $\Gamma\le G$ gives rise to a locally symmetric space $M=\Gamma\setminus X$. A sequence $(M_j)_{j\in\N}$ of locally symmetric spaces $M_j=\Gamma_j\setminus X$ is said to be \textit{Benjamini--Schramm convergent} to $X$ if for every $R>0$, the probability that the ball of radius $R$ centered at a random point in $M_j$ is isometric to a $R$-ball in $X$ tends to $1$ as $j\to\infty.$ A family of lattices in $G$ is \textit{uniformly discrete} if there is a neighborhood of the identity in $G$ that intersects all of its conjugates trivially. 
Benjamini-Schramm convergence has been extensively studied in relation to the asymptotic behaviour of $L^2$-invariants (\cite{abert2017growth}, \cite{abert2011growth}) and quantum ergodicity properties (\cite{le2017quantum}, \cite{brumley2023quantum}, \cite{peterson2023quantum}).
Next, one can associate to any uniform lattice $\Gamma$ in a locally compact group $G$ the spectral measure $\mu_\Gamma$ on the unitary dual $\widehat G$. A sequence of lattices $\{\Gamma_j\}_{j\in\N}$ is said to be \textit{Plancherel convergent} if the sequence $\{\mu_{\Gamma_j}\}_{j\in\N}$ converges to the Plancherel measure of $\widehat G$ in a certain sense, which will be specified in Section \ref{BS and Plancherel}. This kind of convergence has also been studied in many instances (\cite{abert2017growth}, \cite{deitmar2020benjamini}, \cite{wallach1978limit}, \cite{matz2019limit}, \cite{deitmar1999limit}, \cite{finis2015limit}). 
It is known that Plancherel convergence implies Benjamini--Schramm convergence \cite[Thm. 2.6]{deitmar2019benjamini}. Regarding  the converse direction, the authors of \cite{abert2017growth} have shown that a uniformly discrete Benjamini--Schramm sequence is Plancherel convergent. Since uniform discreteness for sequences of lattices in higher rank Lie groups would follow from the Lehmer conjecture, these notions of convergence are expected to be equivalent in the higher rank case. This bears the question whether these two notions of convergence are always equivalent. We will demonstrate that this is not the case by constructing a Benjamini--Schramm convergent sequence, which is not Plancherel convergent. For this we start with the well-known principal congruence subgroups
\begin{equation*}
     \Gamma(N)=\left\{\begin{pmatrix}
           a & b \\ c & d
        \end{pmatrix}\in \SL_2(\Z)\ :\ a,d\equiv 1\ (\text{mod}\ N),\ b,c\equiv 0\ (\text{mod}\ N)
\right\}
\end{equation*}
of $\SL_2(\mathds{Z})$ and replace the cusps of the associated congruence surfaces $X(N)=\Gamma(N)\setminus\mathds{H}$ with short geodesics. By identifying these geodesics in pairs we end up with compact surfaces. An explicit study of the length spectrum of these surfaces allows one to directly check whether they fulfill the above notions of convergence.

\textit{Acknowledgments: }We want to thank Jean Raimbault for the fruitful discussions regarding the above example and sharing his thoughts with us.
Also, we thank Anton Deitmar for reading early drafts of the paper and for helping us through the whole process.

\section{Preliminaries}

In this Section we collect some terminology and known results about locally compact groups and hyperbolic surfaces that we will need in the following.

\subsection{Plancherel and Benjamini--Schramm convergence}\label{BS and Plancherel}

Let $G$ be a connected semisimple Lie group of non-compact type and fix a Haar measure $\mu$ on $G$. Let $\Gamma\le G$ be a uniform lattice and $\widehat G$ be the unitary dual of $G$. Let $\mu_{\Pl}$ be the Plancherel measure on $\widehat{G}$. The right regular representation $R$ on $L^2(\Gamma\setminus G)$ defined by $R_y\phi(x):=\phi(xy)$ decomposes into a direct sum of unitary irreducible representations ,
$$L^2(\Gamma\setminus G)\cong \bigoplus_{\pi\in\widehat G} N_\Gamma(\pi)\pi,$$
with finite multiplicities $N_\Gamma(\pi)$, see e.g. \cite[Thm. 9.2.2]{deitmar2014principles}.
\begin{defn}\label{Plancherel convergence definition}\textbf{(Plancherel convergence): }
    Let $\delta_\pi$ denote the Dirac measure for $\pi\in\hat G$. The measure on $\hat G$ defined by 
    $$
        \mu_\Gamma:=\sum_{\pi\in\widehat G}N_{\Gamma}(\pi)\delta_\pi
    $$
    is called the \textit{spectral measure} associated with $\Gamma$. 
    A sequence $(\Gamma_j)_{j\in\N}$ of uniform lattices in $G$ is called \textit{Plancherel convergent} (or a \textit{Plancherel sequence}) if for every $f\in C^\infty_c(G)$
    \begin{equation}\label{Pl. convergence def 2}
        \frac{1}{\text{vol}(\Gamma_j\setminus G)}\mu_{\Gamma_j}(\hat f)\xrightarrow[]{j\to\infty}\mu_{\Pl}(\hat f).
    \end{equation}
If the sequence $\{\Gamma_j\}_{j\in\N}$ is Plancherel convergent, we also call the associated sequence of locally symmetric spaces $\Gamma_j\setminus G \operatorname{/} K$ Plancherel convergent.   
\end{defn}
Let us turn to Benjamini-Schramm convergence. In the following, $\{\Gamma_j\}_{j\in\N}$ is always a sequence of torsion-free lattices in $G$. 

\begin{defn}
    Let $X=\Gamma\setminus G\operatorname{/} K$ be a locally symmetric space. For $p\in X$, we define the \textit{injectivity radius} at $p$ as 
    $$\InjRad(p)=\sup\Big\{r\in\R\ :\ \operatorname{exp}_p:B_r(0)\to X\ \text{is injective}\Big\},$$
    where $\operatorname{exp}_p:T_p X\to X$ is the exponential map. The \textit{injectivity radius} of $X$ is defined as 
    $$\text{InjRad}(X)=\inf_{p\in X}\operatorname{InjRad}(p).$$
    For $R>0$ we define the \textit{$R$-thin part} of $X$ as 
    $$X_{<R}=\Big\{p\in X\ :\ \text{InjRad}(p)<R\Big\}.$$
\end{defn}
\begin{defn}\label{def of BS-convergence}\textbf{(Benjamini--Schramm convergence): }
    We say that a sequence $(X_j)_{j\in\N}$ of locally symmetric spaces given by $X_j=\Gamma_j\setminus G\operatorname{/}K$ is \textit{Benjamini--Schramm convergent} (or \textit{BS-convergent}) to the universal cover $G\operatorname{/}K$ if for every $R>0$
    \begin{equation}\label{BS-Convergence def.}
        \frac{\text{vol}\left((X_j)_{<R}\right)}{\text{vol}(X_j)}\xrightarrow[]{j\to\infty}0.        
    \end{equation}
    In this case, we will also say that the associated sequence of lattices $(\Gamma_n)_{n\in\N}$ is Benjamini--Schramm convergent.
\end{defn}
Finally, we notice that for $l_{s,j}$ the length of a shortest closed geodesic on $X_j$, also called the \textit{systole} of $X_j$, the sequence $(X_j)_{j\in\N}$ is uniformly discrete if and only if there is a uniform lower bound for the systoles $l_{s,j}$ for all $j\in\N.$

\subsection{Decompositions of hyperbolic surfaces}
In this section we collect known results concerning hyperbolic surfaces and their decomposition into (possibly degenerate) pairs of pants (or $Y$-pieces). For a complete treatment of the topic, we refer to \cite{buser2010geometry}. Let $X$ be a smooth closed hyperbolic surface. Then $X$ is a quotient $X=\Gamma\setminus\Hyp$, where $\Gamma$ is a purely hyperbolic Fuchsian group and $\Hyp$ is the hyperbolic plane 
$$\Hyp=\Big\{z=x+iy\in\C\ :\ y>0\Big\}$$
equipped with the Riemannian metric $ds^2=\frac{dx^2+dy^2}{y^2}$. We write $g=g(X)$ for the \textit{genus} of $X$ and note that Gauss--Bonnet Theorem implies $\text{vol}(X)=4\pi(g-1)$. A compact topological surface has \textit{signature} $(g,n)$ if it is obtained from a closed topological surface of genus $g$ by removing the interior of $n$ disjoint topological disks. A closed hyperbolic surface of signature $(0,3)$ is called a \textit{Y-piece} or a \textit{pair of pants}. For any triple of positive real numbers $l_1, l_2, l_3$ there exists a pair of pants $Y_{l_1,l_2,l_3}$ with boundary geodesics $\gamma_1,\gamma_2,\gamma_3$ of respective lengths $l_1,l_2,l_3$. Such a $Y$-piece is obtained by pasting together two copies of a geodesic hexagon.
\begin{thm}[Collar Theorem] Let $X$ be a smooth closed hyperbolic surface of genus $g$. Let $\gamma_1,\dots, \gamma_m$ be pairwise disjoint simple closed geodesics on $X$. Then:
\begin{enumerate}
    \item $m\le 3g-3$;
    \item There exist simple closed geodesics $\gamma_{m+1},\dots, \gamma_{3g-3}$, which together with $\gamma_1,\dots, \gamma_m$ decompose $X$ into $Y$-pieces;
    \item The collars 
    $$\mathscr{C}(\gamma_i)=\{p\in X\ ;\ \text{dist}(p,\gamma_i)\le w_i\}$$
    of widths 
    $$w_i=\sinh^{-1}\left(\frac{1}{\sinh\left(l_i/2\right)}\right)$$
    are pairwise disjoint for $i=1,\dots, 3g-3$. 
\end{enumerate}
\end{thm}

In order to construct Benjamini--Schramm convergent sequences which are not uniformly discrete, we need to extend our discussion of $Y$-pieces to possibly non-compact surfaces. A \textit{degenerate} hexagon is a hexagon with either one, two or three points at infinity. We will also refer to points at infinity as \textit{punctures}. One can paste together two degenerate hexagons to get a \textit{degenerate $Y$-piece} (or \textit{degenerate pair of pants}). We extend the notation $Y_{l_1,l_2,l_3}$ to degenerate $Y$-pieces by writing $l_i=0$ for any boundary component arising from a puncture. A degenerate $Y$-piece contains a neighborhood $\mathcal{C}$ around each puncture, which is isometric to $(-\infty,\log 2]\times\mathbb{S}^1$ equipped with the Riemannian metric
$$ds^2=dr^2+e^{2r}dt^2.$$
This neighborhood is called a \textit{cusp}. A $Y$-piece has signature $(0,p;q)$ if it has $p$ boundary geodesics and $q$ punctures. A smooth hyperbolic surface of genus $g$ is said to have \textit{signature} $(g,p;q)$ if it has $p$ boundary geodesics and $q$ cusps.

\begin{thm} Let $X$ be a (possibly non-compact) smooth hyperbolic surface of signature $(g,0;q)$. Let $\gamma_1,\dots, \gamma_m$ be pairwise disjoint simple closed geodesics on $X$. Then:
\begin{enumerate}
    \item $m\le 3g-3+q$;
    \item There exist simple closed geodesics $\gamma_{m+1},\dots, \gamma_{3g-3+q}$, which together with $\gamma_1,\dots, \gamma_m$ decompose $X$ into (possibly degenerate) $Y$-pieces;
    \item The collars $\mathscr{C}(\gamma_i)$, $i=1,\dots, 3g-3+q$ and the cusps $\mathcal{C}^1,\dots, \mathcal{C}^q$ are all pairwise disjoint.
\end{enumerate}
\end{thm}

If we pinch a simple closed geodesic $\gamma$ on a smooth hyperbolic surface $X$, i.e. we let $l_\gamma\to0$, then we expect that the collar $\mathscr{C}(\gamma)$ around $\gamma$ converges in a suitable sense to two copies of a cusp $\mathcal{C}$. To make this more precise we need some additional terminology (cf. \cite{buser2014quasiconformal}). Let $Y_{l_1,l_2,l_3}$ be a non-degenerate $Y$-piece with boundary geodesics $\gamma_i,\ i=1,2,3$. For $0\le r_i\le w_i$, The distance sets 
$$\gamma_i^{r_i}=\Big\{ p\in Y\ :\ \operatorname{dist}(p,\gamma_i)=r_i\Big\}$$
are called \textit{equidistant curves}. For degenerate $Y$-pieces the equidistant curves are given by horocycles 
$$h_r=\{p\in \mathcal{C}\ :\ \operatorname{dist}(p,\partial \mathcal{C})=r\}.$$
Now, select in each half-collar or cusp an equidistant curve $\beta_i$ of length $\lambda_i$. Then the closure of the connected component of $Y_{l_1,l_2,l_3}\setminus(\beta_1\cup\beta_2\cup\beta_3)$ not containing any of the boundary geodesics or punctures of $Y_{l_1,l_2,l_3}$ is called a \textit{restricted $Y$-piece} and denoted $Y_{l_1,l_2,l_3}^{\lambda_1,\lambda_2,\lambda_3}$. We will write $Y_{l_1,l_2}^c$ instead of $Y_{l_1,l_2,0}^{l_1,l_2,c}$. A homeomorphism $\phi: Y\to Y'$ of (possibly restricted) $Y$-pieces is called \textit{boundary-coherent} if for corresponding boundary curves $\alpha_i$ of $Y$ and $\alpha_i'$ of $Y'$ in standard parametrization\footnote{The precise description of this terminology can be found in \cite[p.2]{buser2014quasiconformal}}. one has 
$$\phi(\alpha_i(t))=\alpha_i'(t),\qquad \forall t\in[0,1].$$
For each boundary length $l_i$ we let 
$$P_i=\left\{p\in Y_{l_1,l_2,l_3}\ :\ \text{dist}(p,\gamma_i)<\log(2/l_i)\right\}$$
if $0<l_i<2$ and $P_i=\emptyset$ for $l_i\ge 2$. In the degenerate case, we let $P_i$ be the connected component of $Y_{l_1,l_2,l_3}\setminus h_1$ which contains the puncture corresponding to $l_i=0$. Then 
$$\hat Y_{l_1,l_2,l_3}=Y_{l_1,l_2,l_3}\setminus(P_1\cup P_2\cup P_3)$$
is called a \textit{reduced $Y$-piece}. Let us recall that a piecewise smooth map\footnote{A piecewise smooth map is a continuous map which is smooth on the complement of a finite number of curves} ${\Psi:M\to N}$ of Riemannian manifolds $M$ and $N$ is called a \textit{quasi-isometry} if there exists $d>0$ such that for any $p\in M$ and any $v\in T_pM$
\begin{equation}\label{quasi-isometry}
    \frac{1}{d}||v||_M\le||D\Psi(v)||_N\le d||v||_M,
\end{equation}
where $||\cdot||_M$ and $||\cdot||_N$ are the norms associated to the Riemannian metrics on $M$ and $N$ respectively. 
The infimum over all possible $d>0$ such that (\ref{quasi-isometry}) holds is called the \textit{length distortion} and is denoted by $d_\Psi$.

\begin{propx}\label{boundary-coherent homeomorphism} Let $0<l_1,l_2$ and let $0<\eps<\frac{1}{2}$. Set $\eps^*=\frac{2}{\pi}\eps$. Then there exists a boundary-coherent homeomorphism 
$$\phi: Y_{l_1,l_2,\eps}\to Y_{l_1,l_2}^{\eps^*}$$
such that 
\begin{enumerate}
    \item $\phi\left(\hat Y_{l_1,l_2,\eps}\right)=\hat Y_{l_1,l_2,0}$;
    \item The restriction of $\phi$ to $\hat Y_{l_1,l_2,\eps}$ is boundary-coherent and has length distortion $d_\phi\le 1+\frac{5}{4}\eps^2$.
\end{enumerate}
\end{propx}
\begin{proof}
    \cite[Thm. 5.1]{buser2014quasiconformal} 
\end{proof}
\begin{rem}
    The above Proposition can be extended in an obvious manner to $Y$-pieces with more than one degenerating boundary geodesic.
\end{rem}

\section{Criteria for Benjamini--Schramm and Plancherel convergence}\label{section 3}
Given a smooth closed hyperbolic surface $X$ and $R>0$, we denote by $N(X,R)$ the number of closed geodesics on $X$ of length $\le R$ and with $N_s(X,R)$ the number of simple closed geodesics on $X$ with length $\le R$. Let $(X_n)_{n\in\N}$ be a sequence of smooth closed hyperbolic surfaces. By \cite[Prop. 2.2]{raimbault2017convergence}, if 
    \begin{equation}\label{Plancherel condition}
        \frac{N(X_n,R)}{\text{vol}(X_n)}\xrightarrow[]{n\to\infty}0,
    \end{equation}
then $(X_n)_{n\in\N}$ is Benjamini-Schramm convergent. We improve this result by showing that it is enough to consider simple closed geodesics. This will be key in the proof of Theorem \ref{main theorem}. Moreover, we show that (\ref{Plancherel condition}) implies Plancherel convergence.

\begin{propx}\label{criterion for BS convergence}
    Let $(X_n)_{n\in\N}$ be a sequence of closed hyperbolic surfaces. If for every $R>0$
    \begin{equation}\label{Benjamini-Schramm criterion}
        \frac{N_s(X_n,R)}{\operatorname{vol}(X_n)}\xrightarrow[]{n\to\infty}0,
    \end{equation}
    then the sequence $(X_n)_{n\in\N}$ is Benjamini-Schramm convergent.
\end{propx}
\begin{proof} 
The proof of this result is given in section \ref{Appendix 1}. 
\end{proof}

We get the following characterization of Plancherel convergence.

\begin{propx}\label{prop for Pl. convergence}
    A sequence of closed hyperbolic surfaces $(X_n)_{n\in\N}$ is Plancherel convergent if and only if for every $R>0$ 
    \begin{equation}\label{Plancherel characterization}
        \frac{1}{\text{vol}(X_n)}\sum_{\substack{[\gamma]\ne 1 \\ l_\gamma\le R}}\frac{l_{\gamma_0}}{\sinh\left(l_\gamma/2\right)}\xrightarrow[]{n\to\infty}0,
    \end{equation}
    where the sum runs over the length spectrum of $X_n$ and, for each closed geodesic $\gamma$, $\gamma_0$ is the correspondent primitive geodesic.
\end{propx}
\begin{proof}
    The proof of this result is given in section \ref{Appendix 2}. 
\end{proof}

\begin{cor}
    Let $(X_n)_{n\in\N}$ be a sequence of closed hyperbolic surfaces. If for every $ R>0$
$$
        \frac{N(X_n,R)}{\text{vol}(X_n)}\xrightarrow[]{n\to\infty}0,
$$
    then the sequence $(X_n)_{n\in\N}$ is Plancherel convergent.
\end{cor}
\begin{proof}
    Notice that for any closed geodesic $\gamma$ on $X_n$, 
    $$\frac{l_{\gamma_0}}{\sinh\left(\frac{l_\gamma}{2}\right)}\le 2.$$
    Hence,
    $$\frac{1}{\text{vol}(X_n)}\sum_{\substack{[\gamma]\ne 1 \\ l_\gamma\le R}}\frac{l_{\gamma_0}}{\sinh\left(l_\gamma/2\right)}\le 2\frac{N(X_n,R)}{\text{vol}(X_n)}\xrightarrow[]{n\to\infty}0.$$
    Therefore, $(X_n)_{n\in\N}$ is Plancherel convergent. 
\end{proof}

\section{A degenerate Benjamini--Schramm sequence}\label{final section}
This section is to devoted to the proof of

\begin{thm}\label{main theorem}
    There exists a Benjamini--Schramm convergent sequence $(X_n)_{n\in\N}$ of smooth closed hyperbolic surfaces, which is not Plancherel convergent.
\end{thm}
For this, we adapt an example from \cite{buser2006riemann}. Let us recall a few facts about the principal congruence subgroups.
$$
\Gamma(N)=\left\{\begin{pmatrix}
           a & b \\ c & d
        \end{pmatrix}\in \SL_2(\Z)\ :\ a,d\equiv 1\ (\text{mod}\ N),\ b,c\equiv 0\ (\text{mod}\ N)
\right\}.
$$
We denote by $X(N)=\Gamma(N)\setminus\Hyp$ the congruence surface of level $N$. We write $g_N$ and $b_N$ respectively for the genus and the number of boundary components of $X(N)$. The principal congruence subgroup $\Gamma(N)$ is torsion-free for $N\ge3$. All boundary components of $X(N)$ are punctures and we have 
\begin{equation}\label{genus and cusps}
    g_N=1+\frac{d_N(N-6)}{24N},\qquad b_N=\frac{d_N}{2N},
\end{equation}
where $d_N$ is given by $d_2=12$ and $d_N=N^3\prod_{p|N}\left(1-\frac{1}{p^2}\right)$ for $N\ge 3$ (cf. \cite{shimura1971introduction}). Note that the number of cusps of $X(N)$ is always even for $N\ge3$. 
By \cite[Lemma 2]{schmutz1994congruence} the systole $l_{s,N}$ of $X(N)$ is given by 
    \begin{equation}\label{length of systole}
        2\cosh\left(\frac{l_{s,N}}{2}\right)=N^2-2.
    \end{equation}
Now, decompose $X(N)$ into (possibly degenerate) pairs of pants. The boundary components of the pants are either geodesics or punctures. We keep the boundary geodesics and replace each puncture by a geodesic of length $t>0$. Let us reassemble these pieces using the old identifications. Since the number of cusps of $X(N)$ is even, we can identify the remaining geodesics in pairs. This yields a closed hyperbolic surface $X_t(N)$. By counting the number of $Y$-pieces involved one can show that 
\begin{equation}\label{genus of new surface}
    g(X_t(N))\ge g_N.
\end{equation}
We also note that the surface $X_t(N)$ contains $\frac{b_N}{2}$ disjoint simple closed geodesics $\gamma_i,\ i=1,\dots \frac{b_N}{2}$ of length $t$. Now, let $(N_j)_{j\in\N}$ be a sequence of natural numbers such that $N_j\ge3$ and $N_j\xrightarrow[]{j\to\infty}\infty$. Let $(t_j)_{j\in\N}$ be a sequence of positive real numbers converging towards $0$. Let us write $b_j=\frac{b_{N_j}}{2}$ and $g_j=g_{N_j}$. Let $(X_j)_{j\in\N}$ be the sequence of closed hyperbolic surfaces defined by $X_j=X_{t_j}(N_j)$. Also, let $\Gamma_j$ denote the lattice in $G$ such that $X_j=\Gamma_j\setminus \Hyp$.
\begin{propx}\label{Plancherel example}
    The sequence $(X_j)_{j\in\N}$ is Plancherel convergent if and only if $t_j^{-1}$ grows sub-exponentially in $N_j$. 
\end{propx}
\begin{proof}
Let us fix $R>0$ and let $(t_j)_{j\in\N}$ be a sequence converging towards $0$. We claim that for $j$ large enough, any geodesic $\gamma$ in $X_j$ of length $l_\gamma\le R$ is the power of one of the geodesics $\gamma_1,\dots, \gamma_{b_j}$. For a fixed $j\in\N$, let $\gamma$ be a simple closed geodesic in $X_j$, which is not freely homotopic to some power of one of the geodesics $\gamma_1,\dots, \gamma_{b_j}$. 
If $\gamma$ intersects any of the geodesics $\gamma_i,$ $i=1,\dots b_j$ we get from \cite[Cor. 4.1.2]{buser2010geometry} that
    \begin{equation}\label{simple with intersection}
        \sinh\left(l_\gamma/2\right)\ge\frac{1}{\sinh\left(t_j/2\right)}.
    \end{equation}
Hence, $\gamma$ can be dismissed for $j$ large enough. If $\gamma$ does not intersect any of the $\gamma_i,\ i=1,\dots, b_j$, then we have by \cite[Thm. 4.1.1]{buser2010geometry} that $\gamma$ lies outside of the collars $\mathscr{C}(\gamma_i),\ i=1,\dots, b_j$. In this case, there exists a boundary-coherent quasi isometry 
$$\phi:X_j\smallsetminus \bigcup_{i=1}^{b_j} P_i\to X(N_j)\smallsetminus\bigcup_{i=1}^{b_j}P_i'$$
given by the identity on any $Y$-piece where no boundary geodesic has been replaced in the above process, and the map from Theorem \ref{boundary-coherent homeomorphism} in the remaining cases. Its length distortion is bounded by $ d_\phi\le1+\frac{5}{4}t_j^2$. Therefore, $\phi(\gamma)$ defines an element $[\phi(\gamma)]\in\Gamma_j$ such that 
    $$l_{\phi(\gamma)}\le\left(1+\frac{5}{4}t_j^2\right)l_\gamma.$$
    We claim that $[\phi(\gamma)]$ is covered by a hyperbolic transformation, i.e. there exists a closed geodesic in the free homotopy class of $\phi(\gamma)$. Assume otherwise that $[\phi(\gamma)]$ can be covered by a parabolic transformation. Then by \cite[p.72]{abikoff1980real} the curve $\phi(\gamma)$ can be homotoped into the power of a simple loop around a puncture of $X(N_j)$. Now, applying $\phi^{-1}$ gives a homotopy of $\gamma$ into the collar around some geodesic $\gamma_{i_0}$ for $i_0\in\{1,\dots b_j\}$. Hence, $\gamma$ is homotopic to some power of $\gamma_{i_0}$, which is a contradiction to our assumption on $\gamma$. Consequently, there exists a hyperbolic transformation $\eta_\gamma\in\Gamma_j$ which covers $[\phi(\gamma)]$ and we get from (\ref{length of systole}) that
    \begin{equation}\label{bound on simple geodesics}
        \arccosh\left((N_j^2-2)/2\right)\le l_{\eta_{\gamma}}\le l_{\phi(\gamma)}\le\left(1+\tfrac{5}{4}t_j^2\right)l_\gamma.
    \end{equation}
    Since $t_j$ is bounded from above, inequality (\ref{bound on simple geodesics}) shows that for $N_j$ large enough there are no simple closed geodesics $\gamma$ in $X_j$ of length $l_\gamma\le R$ apart from $\gamma_1,\dots,\gamma_{b_j}$.

    Next, let $\gamma$ be a non-simple closed geodesic different from a power of any of the geodesics $\gamma_1,\dots\gamma_{b_j}$. According to \cite[Thm. 4.2.4]{buser2010geometry}, any non-simple primitive geodesic of smallest length is a figure eight geodesic $\delta$ (i.e. a closed geodesic with exactly one self-intersection) embedded into a $Y$-piece. Any $Y$-piece contains at least one boundary geodesic not belonging to $\{\gamma_1,\dots \gamma_{b_j}\}$, as $X(N_j)$ would otherwise not be connected. Then the length formula for $\delta$ (see \cite[eq. (4.2.3)]{buser2010geometry}) yields 
    $$l_\gamma\ge l_\delta\ge 2\left(1+\tfrac{5}{4}t_j^2\right)^{-1}\arccosh\left((N_j^2-2)/2\right).$$
    This proves that for $j$ large enough any geodesic in $X_j$ of length $\le R$ is a power of one of the geodesics $\gamma_1,\dots,\gamma_{b_j}$. 
    Finally, we want to use the characterization of Plancherel convergence given in Proposition \ref{prop for Pl. convergence}. For $j$ large enough, 
    \begin{equation}\label{condition in example}
        \sum_{\substack{[\gamma]\ne 1 \\ l_y\le R}}\frac{l_{\gamma_0}}{\sinh\left(l_\gamma/2\right)}= b_j\sum_{\substack{k\in\N\\ kt_j\le R}}\frac{t_j}{\sinh\left(kt_j/2\right)}=b_j\sum_{\substack{k\in\N\\ kt_j\le R}}\frac{1}{k}\frac{kt_j}{\sinh\left(kt_j/2\right)}\le 2b_j\sum_{\substack{k\in\N\\ kt_j\le R}}\frac{1}{k},
    \end{equation}
    where we used $x\le\sinh(x)$ for $x\in\R$.
    Plugging into (\ref{condition in example}) the asymptotic expansion of the harmonic series
    $$\sum_{k=1}^n\frac{1}{k}=\log n +\gamma_E+\mathcal{O}\left(\frac{1}{n}\right),$$
    where $\gamma_E$ is the Euler-Mascheroni constant, we obtain the estimate 
    \begin{equation}\label{upper bound on pre-Plancherel condition}
        \sum_{\substack{[\gamma]\ne 1 \\ l_y\le R}}\frac{l_{\gamma_0}}{\sinh\left(l_\gamma/2\right)}\le C_1b_j|\log t_j|
    \end{equation}
    for some sufficiently large constant $C_1$. Plugging the values for $b_j$ and $g_j$ from (\ref{genus and cusps}) into (\ref{upper bound on pre-Plancherel condition}) yields
    \begin{equation}\label{upper bound on Plancherel condition}
        \frac{1}{\text{vol}(X_j)}\sum_{\substack{[\gamma]\ne 1 \\ l_y\le R}}\frac{l_{\gamma_0}}{\sinh\left(\frac{l_\gamma}{2}\right)}\le\frac{3C_1|\log t_j|}{\pi(N_j-6)}.
    \end{equation}
   In particular, the right hand side of (\ref{upper bound on Plancherel condition}) converges to $0$ if $t_j^{-1}$ grows sub-exponentially in $N_j$, in which case $(X_j)_{j\in\N}$ is Plancherel convergent.
   Similarly, we get the lower bound
   \begin{equation}\label{lower bound on pre-Plancherel condition}
       \sum_{\substack{[\gamma]\ne 1 \\ l_y\le R}}\frac{l_{\gamma_0}}{\sinh\left(l_\gamma/2\right)}= b_j\sum_{\substack{k\in\N\\ kt_j\le R}}\frac{1}{k}\frac{kt_j}{\sinh\left(kt_j/2\right)}\ge\\
       C_R b_j \sum_{\substack{k\in\N\\ kt_j\le R}}\frac{1}{k}\ge C_Rb_j|\log t_j|,
   \end{equation}
for $C_R$ some positive constant depending on $R$. Hence, we get
    \begin{equation}\label{lower bound on Plancherel condition}
        \frac{1}{\text{vol}(X_j)}\sum_{\substack{[\gamma]\ne 1 \\ l_y\le R}}\frac{l_{\gamma_0}}{\sinh\left(l_\gamma/2\right)}\ge C_R'\frac{|\log t_j|}{N_j-6},
    \end{equation}
    for $C_R'$ only depending on $R$. 
    In particular, the right hand side of (\ref{lower bound on Plancherel condition}) diverges if $t_j^{-1}$ grows at least exponentially in $N_j$. This concludes the argument.
    \end{proof} 

Let us turn to Benjamini--Schramm convergence of the sequence $(X_j)_{j\in\N}$. 

\begin{propx}\label{final prop}
    The sequence $(X_j)_{n\in\N}$ is Benjamini--Schramm convergent if $t_j$ converges to $0$ as $j\to\infty$.
\end{propx} 
\begin{proof}
Let $R>0$ be given. As in the proof of Proposition \ref{Plancherel example} we get that, for $N_j$ large enough, there is no simple closed geodesic in $X_j$ of length $\le R$ except for $\gamma_1,\dots, \gamma_{b_j}$ and they all have length $t_j$. Therefore, we have $N_S(X_j,R)=b_j$ in the terminology of Corollary  \ref{criterion for BS convergence}. Thus, we have for $j$ large enough that
$$\frac{N_S(X_j,R)}{\operatorname{vol}(X_j)}=\frac{b_j}{4\pi(g(X_j)-1)}\le \frac{C_0}{N_j},$$ 
for $C_0$ some positive constant, where we made use equation (\ref{genus and cusps}) and $g(X_j)\ge g_j$. Therefore, 
$(X_j)_{j\in\N}$ is Benjamini--Schramm convergent. 
\end{proof}

The existence of the example from Theorem \ref{main theorem} follows from Proposition \ref{Plancherel example} and Proposition \ref{final prop} by choosing a sequence $(t_j)_{j\in\N}$ such that $t_j^{-1}$ grows at least exponentially in $N_j$.

\section{Appendix} 

\subsection{Proof of Proposition \ref{criterion for BS convergence}}\label{Appendix 1}
For any simple closed geodesic $\gamma$ in $X$ we define 
$$ \mathscr{C}_R(\gamma)=\Big\{ p\in X\ :\ \operatorname{dist}(p,\gamma_0)\le w_R(\gamma_0)\Big\},\quad w_R(\gamma)=\sinh^{-1}\left(\frac{\sinh(R)}{\sinh\left(l_{\gamma}/2\right)}\right).$$
Notice that thanks to the Collar Theorem \cite[Thm. 4.1.1]{buser2010geometry}, for $\tilde R=\sinh^{-1}(1)$ and $w(\gamma)=w_{\tilde R}(\gamma)$, $\mathscr{C}_{\tilde R}(\gamma)$ is isometric to the cylinder $C(\gamma)=[-w(\gamma),w(\gamma_0)]\times \mathbb{S}^1$ with the Riemannian metric $ds^2=d\rho^2+l^2(\gamma)\cosh^2(\rho)dt^2$. In an analogous way we define, for any simple closed geodesic $\gamma$ and $R>0$, the cylinder $C_R(\gamma)=[-w_R(\gamma),w_R(\gamma)]\times \mathbb{S}^1$, endowed with the same Riemannian metric as above. Hence, 
\begin{equation}\label{volume of hyp. cylinder}
    \operatorname{vol}(C_R(\gamma))=4\pi l(\gamma)\sinh(w_R(\gamma))=4\pi \sinh(R)\frac{l_{\gamma}}{\sinh\left(l_{\gamma}/2\right)}.
\end{equation}

\begin{propx}\label{prop for BS convergence}
    Let $X$ be a closed hyperbolic surface and $R>0$. Then 
       $$X_{<R}\subseteq\bigcup_{\gamma\in\mathcal{S}_{2R}(X)}\mathscr{C}_R(\gamma),$$
       where $\mathcal{S}_{R}(X)$ is the set of simple closed geodesics on $X$ with length $\le R$. Moreover, 
    \begin{equation}\label{R-thin part estimate}
        \operatorname{vol}(X_{<R})\le C_R\cdot N_s(X,R),        
    \end{equation}
    where $N_s(X,R)$ is the number of simple closed geodesics on $X$ with length $\le R$ and $C_R=8\pi\sinh(R)$ is a constant only depending on $R$. 
\end{propx}

Estimate (\ref{R-thin part estimate}) on the volume of the $R$-thin part of a hyperbolic surface immediately implies Proposition \ref{criterion for BS convergence}.

\begin{proof}[Proof of Proposition \ref{criterion for BS convergence}]
        By Proposition \ref{prop for BS convergence} we immediately get, for every $R>0$, 
    $$\frac{\operatorname{vol}\left((X_n)_{<R})\right)}{\operatorname{vol}(X_n)}\underbrace{\le}_{(\ref{R-thin part estimate})}C_R\frac{N_s(X_n,R)}{\operatorname{vol}(X_n)}\xrightarrow[]{n\to\infty}0.$$
    Therefore the sequence $(X_n)_{n\in\N}$ is Benjamini--Schramm convergent.
\end{proof}

\begin{figure}[htp]
    \centering
    \includegraphics[width=1\linewidth]{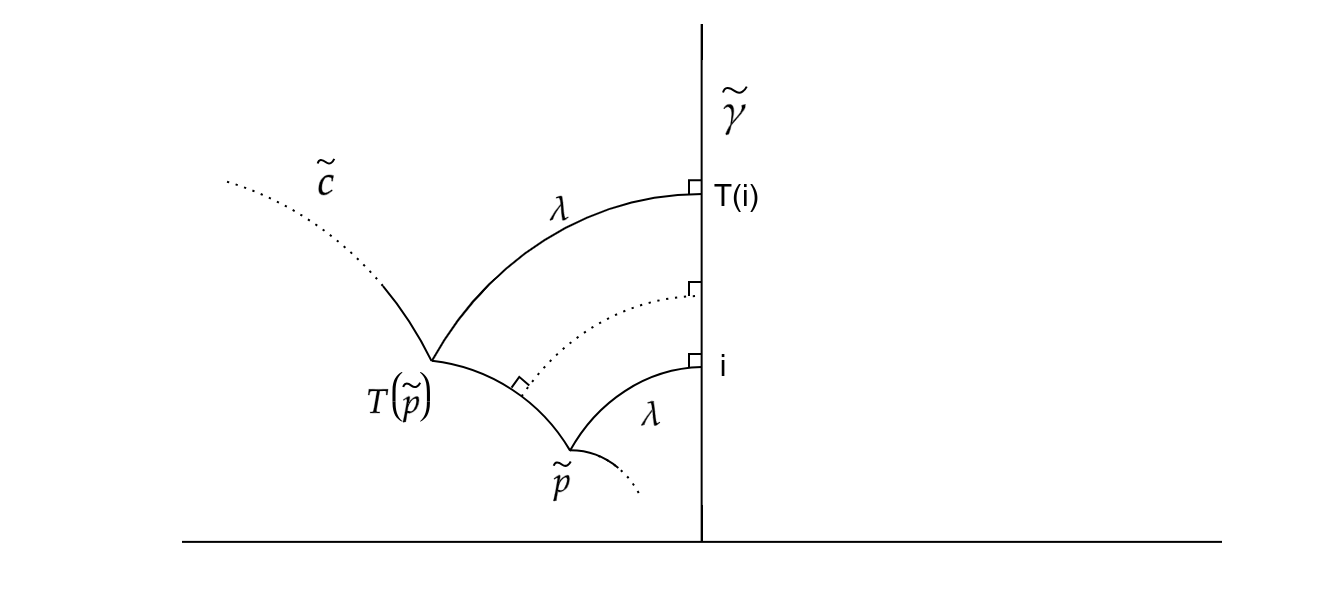}
    \caption{lift to universal cover}
    \label{figure 1}
\end{figure}

We are only left with the proof of Proposition \ref{prop for BS convergence}. For this we will make use of the following result.
\begin{lemma}\label{lemma for BS convergence}
    Let $X$ be a closed hyperbolic surface and $R>0$. If $p\in X$ is such that $InjRad(p)\le R$, then there exists a simple closed geodesic $\gamma$ such that $l_\gamma\le 2R$ and $\operatorname{dist}(p,\gamma)\le w_R(\gamma)$.
\end{lemma}
\begin{proof}
    We adapt the proof of \cite[Thm. 4.1.6]{buser2010geometry}. Let $p\in X$ such that $\text{InjRad}(p)\le R$. There exists a simple geodesic loop $c$ at $p$ of length $2\text{InjRad}(p)$. Let $x\in\gamma$ be the point at minimal distance to $p$, i.e. $\text{dist}(p,\gamma)=\text{dist}(p,x)$. Let $\tilde \gamma$ and $\tilde c$ be homotopic lifts of $\gamma$ and $c$ to $\Hyp$, and let $T$ be the covering transformation with axis $\tilde \gamma$ corresponding to $c$. Without loss of generality, we can assume that $\pi(i)=x$ and $\pi(i\R)=\gamma$. Indeed, for  any lift $\tilde\gamma$ of $\gamma$ and $\tilde x$ in the fiber of $x$ there exists a unique isometry of $\Hyp$ mapping $x$ to $i$ and $\tilde\gamma$ to $i\R$. Hence, $\gamma$ is closed by the covering transformation $T$ associated to $\begin{pmatrix} e^{l_\gamma/2} & 0 \\ 0 & e^{-l_\gamma}/2 \end{pmatrix}\in \SL_2(\R)$. Let $\tilde p\in\tilde c$ be a lift of $p$ and consider its image $T(\tilde p)$. We obtain a geodesic quadrangle with sides of length $\lambda,\ l_\gamma,\ \lambda,\ 2\text{InjRad}(p)$ (see figure \ref{figure 1}). Dropping the common perpendicular between $\tilde c$ and $\tilde \gamma$ we obtain two isometric trirectangles. Standard results from hyperbolic trigonometry (\cite[Thm. 2.3.1]{buser2010geometry}) yield  
    $$\sinh(\text{InjRad}(p))=\sinh(l_\gamma/2)\cosh(\lambda).$$
    Hence, 
    $$\sinh(l_\gamma/2)\le\sinh(\text{InjRad}(p))\le \sinh(R),$$ 
    so that 
    $l_\gamma\le 2R$. Moreover, 
    $$\sinh(l_\gamma/2)\sinh(\lambda)\le \sinh(l_\gamma/2)\cosh(\lambda)\le \sinh(R).$$
    Therefore,
    $$\text{dist}(p,\gamma)\le\lambda\le\sinh^{-1}\left(\frac{\sinh(R)}{\sinh(l_\gamma/2)}\right)=w_R(\gamma).$$
\end{proof}

\begin{proof}[Proof of Proposition \ref{prop for BS convergence}]
    Let $p\in X$ such that $\operatorname{InjRad}(p)\le R$. Then, by \ref{lemma for BS convergence}, there exists a simple closed geodesic such that $\gamma_0\in\mathcal{S}_{2R}(X)$ and $p\in \mathscr{C}_R(\gamma_0)$. Therefore, we obtain 
    $$X_{<R}\subseteq\bigcup_{\gamma_0\in\mathcal{S}_{2R}(X)}\mathscr{C}_R(\gamma_0).$$
    We claim that for any simple closed geodesic $\gamma_0$ on $X$, and for any $R>0$, 
        \begin{equation}\label{volume estimate}
            \operatorname{vol}(\mathscr{C}_R(\gamma_0))\le \operatorname{vol}(C_R(\gamma_0)).
        \end{equation}
    Assuming the Claim, it is immediate to see that 
    $$\operatorname{vol}(X_{<R})\le\sum_{\gamma_0\in\mathcal{S}_{2R}(X)}\operatorname{vol}(\mathscr{C}_R(\gamma_0)\underbrace{\le}_{(\ref{volume estimate})}\sum_{\gamma_0\in\mathcal{S}_{2R}(X)}\operatorname{vol}(C_R(\gamma_0))\underbrace{\le}_{(\ref{volume of hyp. cylinder})}$$
    $$\sum_{\gamma_0\in\mathcal{S}_{2R}(X)}4\pi\sinh(R)\underbrace{\frac{l_{\gamma_0}}{\sinh\left(l_{\gamma_0}/2\right)}}_{\le 2}\le\sum_{\gamma_0\in\mathcal{S}_{2R}(X)}8\pi\sinh(R)=8\pi\sinh(R)N_s(X,2R).$$
    In order to conclude the proof, we are only left with the proof of Claim. Let $\Gamma\le G$ be the uniform lattice in $\SL_2(\R)$ such that $X=\Gamma\setminus \Hyp$. We denote by $\tilde\gamma_0\in\Gamma$ the primitive element of $\Gamma$ corresponding to the geodesic $\gamma_0$ in $X$. Then, $\tilde\gamma_0$ is conjugate to 
        $$\begin{pmatrix}
            e^{l_{\gamma_0}/2} & 0 \\ 0 & e^{-l_{\gamma_0}/2}
        \end{pmatrix}.$$
        Therefore, without loss of generality, we can assume that $\tilde\gamma_0=\begin{pmatrix}
            e^{l_{\gamma_0}/2} & 0 \\ 0 & e^{-l_{\gamma_0}/2}
        \end{pmatrix}$.
        Consider the cylinder $C(\gamma_0)$ given by the quotient $\Gamma_{\tilde\gamma_0}\setminus \Hyp$, where $\Gamma_{\tilde\gamma_0}=\la \tilde\gamma_0\ra$ is the group generated by $\tilde\gamma_0\in SL_2(\R).$
        $C(\gamma_0)$ is a hyperbolic cylinder with only one simple closed geodesic, we denote it $\gamma_0'$ to distinguish it from $\gamma_0$ in $X$, of length $l(\gamma_0')=l(\gamma_0)$. Notice that $\Gamma_{\tilde\gamma_0}$ is a subgroup of $\Gamma$, therefore $C(\gamma)=\Gamma_{\tilde\gamma_0}\setminus\Hyp$ is a covering of $X=\Gamma\setminus \Hyp$. Let $\pi: C(\gamma_0)\to X$ be the covering projection. In particular, we have $\pi(\gamma_0')=\gamma_0$. Furthermore, notice that $C_R(\gamma_0)=[-w_R(\gamma_0),w_R(\gamma_0)]\times\mathbb{S}^1$ is the hyperbolic cylinder obtained by "cutting" $C(\gamma)$ at distance $w_R(\gamma_0)$ from $\gamma_0'$. Let now $p\in \mathscr{C}_R(\gamma_0)$. Then there exists a geodesic arc $\sigma$ joining $\gamma_0$ to $p$, of length $l(\sigma)\le w_R(\gamma_0)$. Since $\pi(\gamma_0')=\gamma_0$, we can lift $\sigma$ to $C(\gamma_0)$ to obtain a geodesic arc $\sigma'$ in $C(\gamma_0)$, joining $\gamma_0'$ to some $p'\in C(\gamma_0)$, of the same length of $\sigma$, i.e. $p'\in C_R(\gamma_0)$. This shows that we can lift $\mathscr{C}_R(\gamma_0)$ via $\pi^{-1}$ into $C_R(\gamma_0)$, hence 
        $$\operatorname{vol}(\mathscr{C}_R(\gamma_0))\le \operatorname{vol} (C_R(\gamma_0)),$$
        which concludes the proof. 
\end{proof}

\subsection{Proof of Proposition \ref{prop for Pl. convergence}}\label{Appendix 2}
In order to prove Proposition \ref{prop for Pl. convergence} we need to introduce some more terminology and recall known results concerning the representation theory of $\SL_2(\R)$. For the following, we refer the reader to \cite[Chapter 11]{deitmar2014principles}. 
The \textit{Hecke algebra} $\mathcal{H}$ of $G=SL_2(\R)$ is the set of $K$-bi-invariant functions on $G$ which are in $L^1(G)$. We call $\mathcal{H}\mathcal{S}$ the subalgebra of the Hecke algebra consisting of smooth functions $f\in\mathcal{H}$ such that $\varphi_f\in\mathcal{S}([0,\infty))$ is a Schwartz function. For every $f\in\mathcal{H}$, there exists a unique function $\varphi_f$ on $[0,\infty)$ such that $f(x)=\varphi_f\left(\operatorname{tr}(x^Tx)-2\right)$ for every $x\in G$. To any $f\in\mathcal{HS}$ we associate a function $h_f(r)=\operatorname{tr}\pi_{ir}(f)$, where the $\pi_{ir}$'s are, for $r\ge0$, the unitary principal series representations of $G$ (\cite{deitmar2014principles}, 11.2.6 + 11.2.9). For $f\in\mathcal{HS}$, one gets the following explicit expression, also known as an explicit Plancherel Theorem (see \cite{deitmar2014principles}, Th. 11.3.1):
\begin{equation}\label{Plancherel measure for SL2}
    \mu_{Pl}(\hat{f})=f(1)=\frac{1}{4\pi}\int_\R h_f(r)\tanh(\pi r)rdr.
\end{equation}
Now let $\Gamma\le G$ be a uniform lattice and let $(r_j)_{j\in\N}$ be a sequence in $\C$ such that $ir_j\in i\R\cup \left(0,\frac{1}{2}\right)$, with the property that $\pi_{ir_j}$ is isomorphic to a subrepresentation of $(R,L^2(\Gamma\setminus G))$ and the value $r=r_j$ is repeated in the sequence $N_\Gamma(\pi_{ir})$ times. For a function $f\in\mathcal{H}$ for which the operator $R(f)$ is of trace class, one gets 
$$\operatorname{tr} R(f)=\sum_{\pi\in\widehat G} N_\Gamma(\pi)\operatorname{tr}\pi(f)=\sum_{j=0}^\infty h_f(r_j).$$
In particular, if $f\in\mathcal{H}\mathcal{S}\cap C^\infty_c(G)$, one can apply Selberg Trace Formula to obtain
\begin{equation}\label{trace formula in hyperbolic geometry}
    \sum_{j=0}^\infty h_f(r_j)=\frac{\operatorname{vol}(\Gamma\setminus \Hyp)}{4\pi}\int_\R h_f(r)\tanh(\pi r)rdr+\sum_{[\gamma]\ne1}\frac{l_{\gamma_0}}{2\sinh\left(l_\gamma/2\right)}g_f(l_\gamma),
\end{equation}
where 
$$g_f(r)=\int_\R \varphi_f\left(2\cosh(r)-2+s^2\right)ds.$$
(\ref{Plancherel measure for SL2}) and (\ref{trace formula in hyperbolic geometry}) yield the following result.
\begin{lemma}\label{Plancherel convergence in SL2}
    Let $(X_n)_{n\in\N}$ be a Plancherel sequence of hyperbolic surfaces. Then, for any $f\in C^\infty_c(G)\cap \mathcal{HS}$, 
    \begin{equation}\label{vanishing sum under Plancherel convergence}
        \frac{1}{\operatorname{vol}(X_n)}\sum_{[\gamma]\ne1}\frac{l_{\gamma_0}}{2\sinh\left(l_\gamma/2\right)}g_f(l_\gamma)\xrightarrow[]{n\to\infty}0,
    \end{equation}
\end{lemma}

We are finally able to prove Proposition \ref{prop for Pl. convergence}.
\begin{proof}[Proof of Proposition \ref{prop for Pl. convergence}]
We start by assuming that (\ref{Plancherel characterization}) holds and prove Plancherel convergence. Let $f\in C^{\infty}_c(G)$ and let us write $K_f:=\operatorname{supp}(f)$. Also, denote 
$$B_n:=\left|\frac{\mu_{\Gamma_n}(\hat f)}{\text{vol}(X_n)}-\mu_{Pl}(\hat f)\right|.$$
Notice that Plancherel convergence of $(X_n)$ is equivalent to $B_n\xrightarrow[]{n\to\infty}0$. Applying Selberg's Trace Formula yields 
$$B_n=\frac{1}{\text{vol}(X_n)}\left|\sum_{[\gamma]\ne 1}l_{\gamma_0}\int_{G_\gamma\setminus G}f(x^{-1}\gamma x)dx\right|.$$

Since the orbital integral is invariant under conjugation, we may assume 
$$\gamma=\begin{pmatrix}
    e^{l_\gamma/2} & 0 \\ 0 & e^{-l_\gamma/2}
\end{pmatrix}.$$
Let $G=ANK$ be the Iwasawa decomposition of $G=SL_2(\R)$. By compactness of $K_f$, there exists $R=R(f)>0$ such that 
$$B_n\le\frac{1}{\text{vol}(X_n)}\sum_{\substack{[\gamma]\ne 1\\ l_\gamma\le R}}l_{\gamma_0}\int_{G_\gamma\setminus G}\left|f\left(x^{-1}\gamma x\right)\right|dx.$$
A quick computation shows that the centralizer $G_\gamma$ of $\gamma$ in $G$ is $G_\gamma=A$, so that 
$$\int_{G_\gamma\setminus G}\left|f\left(x^{-1}\gamma x\right)\right|dx=\frac{1}{2\sinh\left(l_\gamma/2\right)}\int_\R\int_0^{2\pi}\left|f\left(k_\theta^{-1}\begin{pmatrix}
    e^{l_\gamma/2} & y \\ 0 & e^{-l_\gamma/2}
\end{pmatrix}k_\theta\right)\right|d\theta dy.$$
The integrand on the right hand side can only be non-zero if 
$$\begin{pmatrix}
    e^{l_\gamma/2} & y \\ 0 & e^{-l_\gamma/2}
\end{pmatrix}\in KK_fK=:K_f'.$$
$K_f'$ is compact. Hence, there exists $M_f>0$ solely depending on $f$ such that 
$$\int_{G_\gamma\setminus G}\left|f(x^{-1}\gamma x)\right|dx\le \frac{M_f}{\sinh\left(l_\gamma/2\right)}.$$
Thus,
$$B_n\le\frac{M_f}{\text{vol}(X_n)}\sum_{\substack{[\gamma]\ne 1 \\ l_y\le R}}\frac{l_{\gamma_0}}{\sinh\left(l_\gamma/2\right)}.$$
(\ref{Plancherel characterization}) implies that $B_n\xrightarrow[]{n\to\infty}0$, so that $(X_n)_{n\in\N}$ is Plancherel convergent. For the converse, we start by observing that for every $\varphi\in C^\infty_c([0,\infty))$ the function $f(x)=\varphi\left(\operatorname{tr}(x^Tx)-2\right)$ is such that $f\in C^\infty_c(G)\cap\mathcal{HS}$ and $\varphi=\varphi_f$. 

Next, for every $R<0$, choose $\varphi\in C^\infty_c([0,\infty))$ such that $\varphi\ge0$ and $g_f(r)>0$ for every $r\in[0,R]$. Let us denote $m_R:=\min_{[0,R]}g_f(r)$. We get, 
$$\sum_{\substack{[\gamma]\ne 1 \\ l_\gamma\le R}}\frac{l_{\gamma_0}}{\sinh\left(l_\gamma/2\right)}\le\sum_{\substack{[\gamma]\ne 1 \\ l_\gamma\le R}}\frac{l_{\gamma_0}}{\sinh\left(l_\gamma/2\right)}\frac{g_f(l_\gamma)}{m_R}\le \frac{1}{m}\sum_{[\gamma]\ne 1}\frac{l_{\gamma_0}}{\sinh\left(l_\gamma/2\right)}g_f(l_\gamma).$$
By Lemma \ref{Plancherel convergence in SL2}, the last term vanishes as $n\to\infty$, and this concludes the proof. 
\end{proof}

\addcontentsline{}{section}{References}
\bibliographystyle{plain}
\begingroup 

\bibliography{ref} 
\endgroup

\end{document}